\newcommand{\CA}{\hbox{{$\cal A$}}}

\newcommand{\Z}{\mathbb{Z}}
\newcommand{\C}{\mathbb{C}}

\newcommand{\note}[1]{}

\newcommand{\extd}{{\rm d}}

\newcommand{\eps}{{\epsilon}}
\newcommand{\tens}{\mathop{\otimes}}
\newcommand{\id}{{\rm id}}

\renewcommand{\>}{\rangle}

\newcommand{\ad}{{\rm ad}}

\renewcommand{\o}{{}_{\scriptscriptstyle(1)}}
\renewcommand{\t}{{}_{\scriptscriptstyle(2)}}

\newcommand{\proof}{{\bf Proof\ }}
\newcommand{\eproof}{$\quad \diamond$\bigskip}
\newcommand{\eqn}[2]{\begin{equation}#2\label{#1}\end{equation}}

\newcommand{\del}{\partial}
\newcommand{\Dsl}{D\kern-7pt/\ }
\newcommand{\Asl}{A\kern-7pt/\ }
\newcommand{\dsl}{\partial\kern-7pt/\ }

\documentclass[12pt,A4]{article}
\usepackage{amssymb,amsmath}
\newtheorem{lemma}{Lemma}[section]
\newtheorem{propos}[lemma]{Proposition}
\newtheorem{theorem}[lemma]{Theorem}

\begin{document}\baselineskip 18pt

{\ }\qquad \hskip 4.3in \vspace{.2in}

\begin{center} {\Large NONCOMMUTATIVE RICCI CURVATURE
AND DIRAC OPERATOR ON $\C_q[SL_2]$ AT ROOTS OF UNITY}
\\ \baselineskip 13pt{\ }\\
{\ }\\    Shahn Majid\footnote{Royal Society University Research
Fellow}
\\ {\ }\\ School of
Mathematical Sciences\\
Queen Mary, University of London, Mile End Rd\\ London E1 4NS, UK
\end{center}
\begin{center}
12 June -- revised August, 2002
\end{center}

\begin{quote}\baselineskip 14pt
\noindent{\bf Abstract} We find a unique torsion free Riemannian
spin connection for the natural Killing metric on the quantum
group $\C_q[SL_2]$, using a recent frame bundle formulation. We
find that its covariant Ricci curvature is essentially
proportional to the metric (i.e. an Einstein space). We compute
the Dirac operator $\Dsl$ and find for $q$ an odd $r$'th root of
unity that its eigenvalues are given by $q$-integers $[m]_q$ for
$m=0,1,\cdots,r-1$ offset by the constant background curvature. We
fully solve the Dirac equation for $r=3$.

\end{quote}

\baselineskip 17.9pt
\section{Introduction}

This is a follow-up to our work \cite{GomMa:coh} where we covered
the electromagnetic theory on the quantum group $\C_q[SL_2]$ at a
root of unity. We now compute the canonical Riemannian geometry
induced by the Killing form using the frame-bundle formulation of
noncommutative geometry on algebras in
\cite{Ma:rie}\cite{Ma:rieq}. Just as compact Lie groups are
natural examples of Riemannian manifolds, so quantum groups such
as $\C_q[SL_2]$ {\em should} be natural examples of noncommutative
Riemannian geometry. We will see that this is indeed the case if
we use the frame-bundle formulation coming out of quantum group
theory.

There are already a lot of attempts to $q$-deform geometric
structures, including so-called Levi-Civita connections on
specific quantum groups, for example \cite{HecSch:lev}, but on a
largely ad-hoc or case-by-case basis. The problem remains how
definitive or complete any of these various   $q$-formulae might
be. Our goal in the present paper is exactly the opposite: the
frame bundle formulation is defined generally at the level of (in
principle) any unital algebra, including matrix algebras, finite
group function algebras and algebras not connected in any way to
q-deformations. In other words a `functorial theory' at an
algebraic level. Our goal is to take this existing theory and
specialise it to $\C_q[SL_2]$ as a test of the general theory. We
find reasonable answers including a unique spin connection, a
constant Ricci curvature and an interesting Dirac operator $\Dsl$.
We study its modes at roots of unity where all constructions are
finite. The root of unity theory is perhaps the most interesting
both mathematically and physically. It is the case occurring
naturally in the WZW model and hence in some approaches to quantum
gravity\cite{MajSmo:def}. The geometry in this case it little
studied, however; aside from \cite{GomMa:coh}, we note \cite{Coq}
in connection with quantum planes at roots of unity.

We refer to \cite{Ma:rieq} for the frame-bundle formulation but
will recall the key formulae for torsion, metric compatibility,
Ricci tensor etc. when we need them. The main ingredients are a
(possibly noncommutative) algebra $\CA$ of `functions' for the
base space, a quantum group $H$ for the functions on the frame
group, a quantum principal bundle with some algebra $P$ and fixed
subalgebra $\CA$ and of course compatible cotangent bundles
$\Omega^1(\CA),\Omega^1(H),\Omega^1(P)$. Then we need a
representation of the frame quantum group and a framing
isomorphism of $\Omega^1(\CA)$  with the corresponding associated
vector bundle. In the case of interest in the present paper, we
take $\Omega^1(H)=H\tens\Lambda^1$, where $\Lambda^1$ is the dual
of the (braided) Lie algebra of the frame group, and we take
$P=\CA\tens H$ a trivial bundle. The framing isomorphism just
amounts to a left $\CA$-basis $\{e_i\}$ of $\Omega^1(\CA)$ or
`vielbein'. A spin connection just amounts to a Lie-algebra valued
1-form or a map $A:\Lambda^1\to \Omega^1(\CA)$. In our case of a
quantum group $\CA$ we take $H=\CA$ i.e. reduce the frame quantum
group to $\CA$ itself. Then there is a natural choice of vielbein
as explained in \cite{Ma:rieq}, given by the Maurer-Cartan form.
So the new result of the present paper is to obtain a unique
compatible spin connection and thereby compute the entire
Riemannian geometry at least for $\C_q[SL_2]$. We expect similar
results for other standard $\C_q[G]$. Let us mention one
interesting phenomenon; as for the finite group $S_3$ in
\cite{Ma:rieq} we find that the Ricci tensor is essentially
proportional to the Killing metric (an `Einstein space') but
shifted by an invariant element $\theta\tens\theta$. Here $\theta$
is a `time-direction' induced by the noncommutative geometry and
not part of the 3-geometry of $\C_q[SL_2]$. The shift by
$\theta\tens\theta$ turns out to be exactly such that the Ricci
tensor has zero modes.

It is also a very interesting question how such examples as
$\C_q[SL_2]$ could fit into the mathematically more sophisticated
K-theory and operator algebra approach of A. Connes and
others\cite{Con}. There the geometry is defined by a spectral
triple or `axiomatic Dirac operator' obeying some axioms suggested
by operator theory, including self-adjointness. By contrast we
find that our constructed Dirac operator on $\C_q[SL_2]$ at an odd
root of unity has eigenvalues on a complex circle with properties
perhaps reminiscent of the discrete series Jones index of type
$II_1$ subfactors.

\section{Preliminaries}

Here we fix notations in the conventions that we will use and do
some preliminary computations, in Section~2.2. We let $q^2\ne 1$.
The quantum group $\CA=\C_q[SL_2]$ has a
matrix of generators $\{t^\alpha{}_\beta\}=\begin{pmatrix}a& b\\
c& d\end{pmatrix}$ with relations
\[ ba=qab,\quad ca=qac,\quad db=qbd,\quad dc=qcd,\quad cb=bc,
\quad da-ad=q\mu bc,\quad ad-q^{-1}bc=1,\] where $\mu=1-q^{-2}$.
The coproduct $\Delta$ and counit $\eps$ have the usual matrix
coalgebra form. We denote the antipode or `matrix inverse' by $S$.

We will also work with the dimensional Hopf algebra
$\CA=\C_q[SL_2]$ reduced at $q$ a primitive $r$'th root of unity.
This has the further relations
\[ c^r=b^r=0,\quad a^r=d^r=1\]
where $a^r,b^r,c^r,d^r$ generate an undeformed $\C[SL_2]$ central
sub-Hopf algebra of the original $\C_q[SL_2]$. Note also that in
the reduced case $a=(1+q^{-1}bc)d^{-1}$ is redundant and
$\dim(\CA)=r^3$. In this case a basis of $\CA$ is $\{c^kb^nd^m\}$
for $0\le m,n,k\le r-1$.  Explicit  computations are done via
Mathematica for $r=3,5,7$ for concreteness, but we expect
identical results for all odd $r$. Equations involving only the
invariant differential forms do not directly involve the function
algebra and are solved for all $q^2\ne -1$.

\subsection{Exterior algebra}

We take the standard bicovariant exterior algebra on $\C_q[SL_2]$
which has the lowest dimensional (4d) space of
1-forms\cite{Wor:dif}. This is the same set-up as in
\cite{GomMa:coh} except that we carefully switch to $\Lambda^1$
the space of left-invariant 1-forms as used in \cite{Ma:rieq}
(rather than right-invariant). Thus, we take a basis
$\{e_\alpha{}^\beta\}=\begin{pmatrix}e_a& e_b\\
e_c & e_d\end{pmatrix}$, where $e_1{}^2=e_b$, etc., and form a
right crossed module with right multiplication and the right
coaction
\[ \Delta_R(e_\alpha{}^\beta)=e_\gamma{}^\delta\tens
t^\gamma{}_\alpha St^\beta{}_\delta.\] $\Omega^1=\CA\tens
\Lambda^1$ is spanned by the left-invariant forms as a free left
module over $\C_q[SL_2]$. These also generate the invariant
exterior algebra $\Lambda$ and $\Omega=\CA\tens\Lambda$. Here
$e_a,e_b,e_c$ behave like usual forms or Grassmann variables and
\[ e_a\wedge e_d+e_d\wedge e_a+\mu  e_c\wedge e_b=0,\quad
e_d\wedge e_c+q^2e_c\wedge e_d+\mu  e_a\wedge e_c=0\] \eqn{wedge}{
e_b\wedge e_d+q^2e_d\wedge e_b+\mu  e_b\wedge e_a=0,\quad
e_d^2=\mu e_c\wedge e_b.} It is useful to define  $e_z\equiv
q^{-2}e_a-e_d$, which obeys \eqn{wedgez}{ e_z\wedge e_c+q^2
e_c\wedge e_z=0,\quad e_b\wedge e_z+q^2 e_z\wedge e_b=0,\quad
e_z\wedge e_z=(1-q^{-4})e_c\wedge e_b.}
 The relations among 1-forms are obtained by setting to zero the
 kernel of $\Psi-\id$, where $\Psi$ is the crossed-module braiding. Its specific form
can be obtained by converting the right-invariant formulae in
\cite{GomMa:coh} to left ones by $e_a\leftrightarrow e_d$ and
left-right reversal of all products including tensor products. We
will need the eigenvectors of $\Psi$ in Section~4.

The right module structure on 1-forms is defined via the
commutation relations
\[ e_a
\begin{pmatrix}a&b\\ c&d\end{pmatrix}=\begin{pmatrix}qa&q^{-1} b\\
qc&q^{-1}d\end{pmatrix}e_a,\quad [e_c, b]=[e_c,
d]=[e_b,a]=[e_b,c]=0\]
\[[ e_c, a]=q\mu \, be_a,\quad [e_c, c]=q\mu \, de_a,\quad
[e_b, b]=q\mu \, ae_a,\quad [e_b, d]=q\mu \, ce_a\]
\[ [e_d,a]_{q^{-1}}=\mu  be_b,\quad
[e_d,b]_q=\mu  ae_c+q\mu ^2be_a,\quad [e_d,c]_{q^{-1}}=\mu
de_b,\quad [e_d,d]_q=\mu  ce_c+q\mu^2d e_a\] where
$[x,y]_q=xy-qyx$   and $a,b,c,d\in \C_q[SL_2]$.

\subsection{Exterior derivative and Lie bracket structure constants}

The differential on the exterior algebra structure is defined by
graded anticommutator $\extd =\mu^{-1}[\theta,\ \}$ where
$\theta=e_a+e_d$. In particular \eqn{extd1}{ \extd e_a=-e_c\wedge
e_b,\quad \extd e_d=e_c\wedge e_b,\quad \extd e_c=q^2e_c\wedge
e_z,\quad \extd e_b=-e_b\wedge e_z,\quad \extd
e_z=-(q^{-2}+1)e_c\wedge e_b.}

\begin{lemma}cf.\cite{GomMa:coh} For all invertible $q^2\ne 1$,
\begin{eqnarray*} e_a.c^kb^nd^m &=& q^{-m-n+k}c^kb^nd^m.e_a\\
e_b.c^kb^nd^m &=&
c^kb^nd^m.e_b+q^2\mu[n]q^{-m-n}c^kb^{n-1}d^{m-1}e_a+q^{1-n}\mu[m+n]c^{k+1}b^nd^{m-1}e_a\\
 e_c.c^kb^nd^m&=& c^kb^nd^m.e_c+q^{k-m}\mu [k]c^{k-1}b^nd^{m+1}e_a\\
 e_d.c^kb^nd^m&=& q\mu^2([k+1][m+n]c^kb^nd^m+q^{-m}[k][n]c^{k-1}b^{n-1}d^m).e_a+\mu [k] q^n c^{k-1} b^n
d^{m+1}e_b\\
&&+q^{-k}\mu\left([n]c^kb^{n-1}d^{m-1}+q^{-1}(q^{m+n}[m]
+[n])c^{k+1}b^nd^{m-1}\right)e_c+q^{-k+n+m}c^kb^nd^m.e_d
\end{eqnarray*}
where terms with negative powers of $b,c$ are omitted.
 \end{lemma}
\proof  This is a left-handed version of \cite{GomMa:coh}, so we
will
be brief. Writing  $[n]=(q^n-q^{-n})/(q-q^{-1})$ and $u=\begin{pmatrix}a\\
c\end{pmatrix}$, $v=\begin{pmatrix}b\\ d\end{pmatrix}$, we iterate
the bimodule relations to obtain
\[\begin{array}{rll}e_d v^n &=\  q^{n}\,  v^ne_d+q^{n-1}\, [n]\,
\mu\, u v^{n-1} e_c
    +q\mu^2\, [n]\, v^n e_a ,\quad  & e_b v^n \ =\  v^n e_b+q\mu\, [n]\, uv^{n-1} e_a \\
 e_c v^n&=\ v^n e_c ,  & e_a v^n\ =\  q^{-n}v^n e_a \\
 e_d u^n&=\ q^{-n}\, u^n e_d +q^{1-n}\, [n]\, \mu\, v u^{n-1} e_b, &  e_b u^n\ =\   u^ne_b
\\
 e_cu^n &=\ u^n e_c +q\mu\, [n]\, v u^{n-1} e_a  , &e_a u^n \
=\ q^{n}\, u^n e_a.\end{array}\]These then give the commutation
relations with basis elements as stated. \eproof

From these it is clear how the right handed derivatives in
\cite{GomMa:coh} transpose to our left handed ones. We easily
obtain
\begin{eqnarray}\label{genextd}
\extd(c^k b^n d^m)&=&  \quad \mu^{-1}(q^{m+n-k}-1)\, c^k b^n d^m
e_d+
 q^{n-k+1}[k]_{q^2}\, c^{k-1} b^n\, d^{m+1}\, e_b \nonumber
\\
    && + q^{-k-n}\,
(\, [m+n]_{q^2}\, c^{k+1} b^nd^{m-1}+q[n]_{q^2}\,
c^kb^{n-1}d^{m-1})e_c   \nonumber
\\
&& +\mu\, q^{-k-m-n+2}(\, [k+1]_{q^2}\, [m+n]_{q^2}\,  c^k b^n d^m
+
        q[n]_{q^2}\, [k]_{q^2}\, c^{k-1}b^{n-1}d^m) e_a
        \nonumber
\\
    && +\mu^{-1}(q^{-m-n+k}-1)\, c^k b^n  d^m e_a
\end{eqnarray}
obeying the Leibniz rule. We write $[n]_{q^2}=(1-q^{2n})/(1-q^2)$.

Next, we will need the projection $\tilde\pi:\CA\to \Lambda^1$
which characterises the above calculus as a quotient of the
universal one, i.e. with $\extd f=f\o \tilde\pi(f\t)$ for all
$f\in \CA$. Here $\Delta f=f\o\tens f\t$ is the Sweedler notation
for the coproduct. Actually $\tilde\pi$  can be obtained backwards
from $\extd$ as follows. Let the partial derivatives
$\del^i:\CA\to\CA$ be defined by $\extd f=\sum_i (\del^i f)e_i$.
Then \eqn{tildepi}{ \tilde\pi(f)=\sum_i e_i\eps(\del^i f).} In
particular, we obtain \eqn{tildepia}{
\tilde\pi(a)={q\over[2]_q}(qe_a-e_d),\quad \tilde\pi(c)=e_b,\quad
\tilde\pi(b)=e_c,\quad \tilde\pi(d)={1\over[2]_q}(q^2
e_d+(q^2-q^{-1}-1)e_a).} We use the formula (\ref{genextd}) for
the exterior derivative. For generic $q$ we compute $\extd a$
separately.

Finally, we need the braided-Lie algebra structure constants
\eqn{adLR}{ \ad_R=(\id\tens\tilde\pi)\Delta_R,\quad
\ad_L=(\tilde\pi\tens\id)\Delta_L} where
$\Delta_L=(S^{-1}\tens\id)\circ\tau\circ\Delta_R$ is the right
coaction converted to a left coaction and $\tau$ is the usual
vector space flip. We have
\begin{eqnarray}\label{ad} \ad_R(e_a)&=&e_c\tens e_b-e_b\tens e_c
+\mu (e_a-q^2 e_d)\tens
e_a\nonumber\\ \ad_R(e_d)&=&e_b\tens e_c-e_c\tens e_b-\mu(e_a-q^2
e_d)\tens e_a \nonumber\\
\ad_R(e_b)&=&e_b\tens(q^2e_a-e_d)-
(q^{-2}e_a-e_d)\tens e_b \nonumber\\
\ad_R(e_c)&=&e_c\tens ((q^2-1-q^{-2})e_a+q^2 e_d)+
(e_a-q^2 e_d)\tens e_c \nonumber\\
 \ad_L(e_a)&=&e_c\tens e_b-e_b\tens e_c+\mu e_a\tens (e_a-q^2 e_d)
 \nonumber\\
\ad_L(e_d)&=&e_b\tens e_c-e_c\tens e_b-\mu e_a\tens (e_a-q^2e_d)
\nonumber\\
 \ad_L(e_c)&=&(q^2e_a-e_d)\tens e_c-
e_c\tens
(q^{-2}e_a-e_d)\nonumber\\
\ad_L(e_b)&=&((q^2-1-q^{-2})e_a+q^2 e_d)\tens e_b+ e_b\tens
(e_a-q^2 e_d)\end{eqnarray} We define by
$\ad(e_i)=\sum_{j,k}\ad(jk|i)e_j\tens e_k$, where we use the
indices $i,j,k$ to run through $\{e_a,e_b,e_c,e_d\}$. They are
given as above for all $q$ and are versions of the structure
constants of the braided-Lie algebra $\widetilde{sl_{q,2}}$ as
explained in \cite{Ma:rieq}.

\section{Spin connection and Riemannian curvature}

Still with $q$ arbitrary, there is a natural $\Delta_R$-covariant
Killing metric in \cite{Ma:rieq} of the form
\[\eta=e_c \tens e_b + q^2
e_b\tens e_c + {(e_a \tens e_a - q e_a \tens e_d - q e_d \tens e_a
+ q(q^2 + q - 1)e_d \tens e_d)\over
[2]_q}+\lambda\theta\tens\theta\] where $\lambda=q(1- q - q^2)/(1
+ q)$ is the natural choice for the Hodge * operator as explained
in \cite{GomMa:coh}.  In this case we have more simply
\eqn{metric}{ \eta=e_c\tens e_b+q^2 e_b\tens e_c
+{q^4\over[2]_{q^2}}(e_z\tens e_z-\theta\tens\theta)} which (for
real $q$) is the usual $q$-Minkowski space metric cf.
\cite{Ma:book}. This is the local cotangent space of $\C_q[SL_2]$
with $\theta$ an intrinsic `time' direction induced by
noncommutative geometry. We can add any multiple of
$\theta\tens\theta$ and still retain $\Delta_R$-invariance.

For any such invariant metric, we have symmetry in the sense
\eqn{symm}{ \wedge(\eta)=0.} Moreover, the equations for a
torsion-free and skew-metric-compatible `generalised Levi-Civita'
 spin connection become independent of $\eta$ and reduce to the
torsion and `cotorsion' equations\cite{Ma:rieq}: \eqn{torcotor}{
\extd e_i+ \sum_{j,k}A_j\wedge e_k\, \ad_L(jk|i)=0,\quad \extd
e_i+\sum_{j,k}e_j\wedge A_k\, \ad_R(jk|i)=0,\quad \forall i.} In
these equations we write $A(e_i)\equiv A_i$, and a generalised
spin connection is given by four such 1-forms $A_a,A_b,A_c,A_d$
obeying (\ref{torcotor}). In principle there is also an optional
`regularity' condition as explained in \cite{Ma:rieq} which
ensures that the curvature is braided-Lie algebra valued. By the
same arguments as above, this regularity condition can be written
as \eqn{reg}{ \sum_{i,j}A_i\wedge A_j\, \eps(\del^i\del^j
f)=0,\quad \forall f\in \ker\tilde\pi.}

\begin{theorem} For generic $q$ or for $q$ an odd root of unity
there is a unique torsion-free and cotorsion-free spin connection
given by
\[ A_a={q^4\over [4]_{q^2}}e_z=-A_d,\quad A_b
={1\over q^2+q^{-2}}e_b,\quad A_c=
 {1\over q^2+q^{-2}}e_c.\]
The connection is not in general regular.
\end{theorem}
\proof We define the components of the connection by $A_i=A^j{}_i
e_j$ and have to solve for this $4\times 4$ matrix of elements of
$\C_q[SL_2]$. Looking first at the torsion equation, we see that
the coefficients of $A$ are all to the left and hence its
functional dependence is immaterial. We write out the equations
using the form of $\ad_L$ and $\extd e_a$ from (\ref{ad}) and
match coefficients of a basis of $\Lambda^2$. This is a linear
system and we obtain the general solution of the torsion equation
\[ A=\begin{pmatrix}x & q^2\mu y&-q^{-2}\mu z&-q^{-2}\delta\\
z& { q^4 +x (1-q^8)\over 1 + q^4}& 0&  q^{-6}z\\
 y&0& {1 +
x(q^6-q^{-2})\over 1+q^4}&
(q^4-1+q^{-2})y\\
-q^2x& -q^4\mu y& \mu z& \delta\end{pmatrix};\quad
\delta={1+q^6\over[4]_{q^2}}- x(q^4-q^2+q^{-2})\] where $x,y,z$
are arbitrary elements of $\C_q[SL_2]$. If we assume for the
moment that $A^j{}_i$ are numbers (i.e. commute with the 1-forms)
then the cotorsion equation is entirely similar, now with $\ad_R$.
Inserting the above gives a joint solution of this form iff
$y=z=0$ and $x=q^2/[4]_{q^2}$ which then simplifies to the
solution stated in the theorem. We still have to show uniqueness.
Denoting our stated solution by $A_0$, our above general solution
of the torsion equation can be written as $A=A_0+\tau$, where
\begin{eqnarray*} \tau_a&=&q^2xe_z+z e_c +y e_b,\quad
\tau_d=x(q^4-q^2+q^{-2})e_z+q^{-6}ze_c+(q^4-1+q^{-2})ye_b\\
\tau_b&=&-\mu z e_z-q^{-2}x(1-q^4)e_b,\quad \tau_c=q^4 \mu y
e_z+x(1-q^4)e_c.\end{eqnarray*} Here we have relabelled
$x-q^2/[4]_{q^2}\to x$ and $x,y,z$ are arbitrary elements of
$\C_q[SL_2]$. We now write out the full cotorsion equation for
$\tau$, which refers only to the linear term since $A_0$ is
already a solution. From $\ad_R(e_a)$, $\ad_R(e_c)$ and
$\ad_R(e_b)$ we have respectively the equations
\begin{eqnarray*}
0&=&-e_c\wedge (\mu z e_z+q^{-2}(1-q^4)xe_b)-e_b\wedge(q^4\mu y
e_z+x(1-q^4)e_c)+\mu q^2 e_z\wedge(q^2x e_z+z e_c+y e_b)\\
0&=&(q^4-1)e_c\wedge ( q^2xe_z+q^{-2}\mu
ze_c+(q^2+q^{-2})ye_b)+q^2 e_z\wedge(q^4 \mu y
e_z+x(1-q^4)e_c)\\
0&=& (q^4-1)e_b\wedge (q^{-2}xe_z+q^{-4}(q^2+q^{-2})ze_c
-\mu y e_b)+e_z\wedge(\mu z e_z+q^{-2}(1-q^4)xe_b).
\end{eqnarray*} and a similar equation from
 We write $e_z f\equiv \sum_i \Delta_z{}^i(f)e_i$ for
any element $f\in \C_q[SL_2]$, and similarly for $\Delta_b{}^i$
and $\Delta_c{}^i$. Using this notation we pass all functions to
the left and then collect basic two forms in a basis $\{e_{ab},
e_{ac}, e_{bd}, e_{cd}, e_{ad}, e_{bc}\}$, where $e_{ab}\equiv
e_a\wedge e_b$ is a shorthand. We use the relations from
Section~2.1 to order products of 1-forms. This gives us up to 18
equations for $x,y,z$. Among these, the coefficient of the
$e_{bd}$ in the first displayed equation above and that of
$e_{bc}$ in the second, and a similar combination from the first
and third equations, gives \eqn{gradyz}{
\Delta_z{}^d(y)={1+2q^4\over 1-q^2+q^4}\, y,\quad
\Delta_z{}^d(z)={2+q^4\over 1-q^2+q^4} z.} But from Lemma~2.1 and
a similar formula for powers of $a$, we have
$-\Delta_z{}^d(a^lc^kb^nd^m)=q^{m+n-k-l}a^lc^kb^nd^m$. Here
$m+n-k-l$ is the natural $\Z$-grading on $\C_q[SL_2]$,
corresponding to the diagonal $\C^*$ or $U(1)$ action. Taking a
homogeneous basis of $\C_q[SL_2]$ (such  as in Lemma~2.1 along
with  another part with $a$ to the left in place of $d$ on the
right for the full non-reduced Hopf algebra), and noting that for
typical $q$ the equation $-(1+2q^4)=(1-q^2+q^4)q^N$ has no
solution among integers $N$, we conclude that $y=0$. Similarly
$z=0$. Also from our 18 equations we have \eqn{gradx}{
\Delta_z{}^b(z)=(1+q^{-2})(\Delta_z{}^d(x)-q^2(x)),\quad
\Delta_z{}^c(y)=(1+q^2)(\Delta_z{}^d(x)-q^{-2}x)} from which we
conclude $x=0$ as well. These arguments apply at least for real
positive $q$ and for $|q|=1$ other than an eighth-root of unity,
since then the eigenvalue equations (\ref{gradyz}) cannot be
solved by integer gradings (by looking at the sign or the modulus
of the two sides). This covers the cases of particular interest.
For all generic $q$ with $|q|\ne 1$ we proceed as follows. (i) If
$y=0$ or $z=0$ then (\ref{gradx}) implies that $x$ is an
eigenfunction of the degree operator, which is not possible as
$q^P+q^{\pm 2}=1$ has no solution among integers $P$, hence $x=0$.
Two more of our 18 equations are \eqn{gradyzx}{
\Delta_z{}^b(x)=\Delta_z{}^d(y){(1+q^2)(-1+q^2-2q^4+q^6)\over
q^4(1+2q^4)},\quad
\Delta_z{}^c(x)=-\Delta_z{}^d(z){(1+q^2)(-1+2q^2-q^4+q^6)\over
q^4(2+q^4)}} and these now imply all $x,y,z=0$. (ii) It remains to
consider both $y,z\ne 0$, so that (\ref{gradyz}) holds with $N,M$
the degrees of $y,z$. Then (\ref{gradyzx}) implies $x\ne 0$ also.
From Lemma~2.1 we see that $\Delta_z{}^b$ raises degree by 2, and
$\Delta_z{}^c$ lowers by two. One may then argue from
(\ref{gradx}) that $x$ is homogeneous of degree $P=M+2=N-2$. Hence
$N-M-4=0$ and the ratio of the eigenvalue equations (\ref{gradyz})
becomes $2+q^{-4}=2+q^4$, which contradicts our assumption that
$|q|\ne 1$. Hence $x=y=z=0$, completing our proof of uniqueness.
The odd roots of unity case was already covered in the above
analysis. It can also be verified
 at 3,5,7'th roots of unity by encoding the relations in
Lemma~2.1 as a $4r^3\times 4r^3$ matrix and computing directly.
The failure of regularity is also verified at a small root of
unity but is expected for all $q\ne 1$. \eproof

The covariant derivative $\Omega^1\to \Omega^1\tens_{\CA}\Omega^1$
is computed from\cite{Ma:rieq}
 \eqn{nabla}{ \nabla e_i=-\sum_{j,k} A_j\tens
e_k\, \ad_L(jk|i),\quad \forall i.} It obeys the usual
derivation-like rule for covariant derivatives, so we need only
give it on basic 1-forms. For the above canonical spin connection
it comes out (in a similar manner to solving the torsion equation
in Theorem~3.1) as
\begin{eqnarray}\label{exnabla} \nabla e_a&=&-\nabla e_d
={1\over q^2+q^{-2}}(e_b\tens
e_c-e_c\tens e_b)- \mu{q^6\over[4]_{q^2}}
e_z\tens e_z\nonumber \\
 \nabla e_b&=& {q^4(1+q^{-2})\over[4]_{q^2}}e_z\tens e_b-{q^2\over
q^2+q^{-2}} e_b\tens e_z\nonumber\\
 \nabla e_c&=&-
{q^4(1+q^{2})\over[4]_{q^2}} e_z\tens e_c+{1\over q^2+q^{-2}}
e_c\tens e_z.\end{eqnarray}

Note that $\nabla \theta=0$, i.e this form is covariantly
constant. It is easy to verify that $\nabla\wedge=\extd$ on
1-forms (i.e. the torsion indeed vanishes). One may likewise also
verify that the cotorsion indeed vanishes, i.e. that skew-metric
compatibility holds in the sense\cite{Ma:rie} \eqn{skewmet}{
(\nabla\wedge\id-\id\wedge\nabla)\eta=0} for $\eta$ of the form
above. On the other hand one may compute that if we extend
$\nabla$ to tensor powers as a derivation (while keeping the
left-most output of $\nabla$ to the far left), then
\eqn{nablacomp}{ \nabla \eta=O(\mu)} i.e. for generic $q$ it tends
to 0 as $q\to 1$, but is not zero for $q^2\ne 1$. This confirms
the view in \cite{Ma:rie} that naive metric compatibility
$\nabla\eta=0$ is too strong for noncommutative geometry and the
correct notion is the skew version (\ref{skewmet}) as expressed in
the cotorsion. We see that there is still a unique solution and it
obeys $\nabla \eta=0$ in the classical limit. Indeed, $A_z\equiv
q^{-2}A_a-A_d,A_b,A_c$ in Theorem~3.1 tends to the classical spin
connection on $SU_2$ when suitable linear combinations of
$A_b,A_c$ are taken. As explained in \cite{Ma:rie},  on the
compact real form of a Lie group $G$ associated to a simple Lie
algebra $g$ it is enough to take the framing group as $G$ itself
and $\Lambda^1=g^*$, with a canonical spin connection given by
${1\over 2}e$, where $e$ is the Maurer-Cartan form; the invariant
Killing form means that ${\rm Ad}(G)\subseteq SO(\dim(g))$ and
$\nabla$ from the connection on this reduced frame bundle is the
usual Levi-Civita metric-compatible one. The unique solution in
Theorem~3.1 is exactly a $q$-deformation of this classical
construction.

We also note from the theorem that if one wants regularity then
one must give up either the torsion or cotorsion conditions, or
both. In fact it is the difference between the torsion and
cotorsion which is the physically important skew
metric-compatibility (\ref{skewmet}), as explained in
\cite{Ma:rieq}, i.e. we should at least keep the torsion and
cotorsion equal. Among connections with constant components there
is a 6-dimensional space of these, of which $A_b\propto e_c$
(other components zero), and a similar one with $A_c\propto e_b$,
are the only regular ones. The moduli of all connections with
torsion=cotorsion is much bigger (e.g. for $r=3$ it is
74-dimensional), among which one may expect some interesting but
not constant regular ones.

Finally, for any connection, the Riemannian curvature is computed
from\cite{Ma:rieq}
 \eqn{riem}{ {\rm
Riemann}=((\id\wedge\nabla)-\extd\tens\id)\circ\nabla} or equally
well from the curvature $F=\extd A+A*A:\ker\eps\to \Omega^2$ of
$A$. When the connection is not regular, the latter Yang-Mills
curvature does not descend to a map $\Lambda^1\to\Omega^2$ (it is
not `Lie algebra valued but lives in the enveloping algebra of the
braided-Lie algebra). However, this does not directly affect the
Riemannian geometry (it merely complicates the geometry `upstairs'
on the quantum frame bundle); in the proof of
\cite[Corol.~3.8]{Ma:rieq} one should simply omit the $\tilde\pi$
in the argument of $F$ for the relation to the Riemann curvature.

\begin{propos} The Riemann curvature $\Omega^1\to
\Omega^2\tens_{\CA}\Omega^1$ of the canonical spin connection in
Theorem~3.1 is
\begin{eqnarray*} {\rm Riemann}(e_a)&=&-{\rm Riemann}(e_d)={q^6\over
(1+q^4)^2}(e_c\wedge e_z\tens e_b+q^{-2}e_b\wedge e_z\tens e_c-\mu
e_c\wedge e_b\tens e_z)\\
 {\rm Riemann}(e_b)&=&{q^6\over (1+q^4)^2}(-e_b\wedge e_z\tens
e_z+q^{-2}(1+q^{-2})e_c\wedge e_b\tens e_b)\\
{\rm Riemann}(e_c)&=&{q^6\over (1+q^4)^2}(-e_c\wedge e_z\tens
e_z-(1+q^{-2})e_c\wedge e_b\tens e_c)\end{eqnarray*}
\end{propos}
\proof Direct computation using the relations in the preliminaries
and the formula (\ref{exnabla}) for $\nabla$. Note that Riemann is
a tensor, so that ${\rm Riemann}(fe_a)=f{\rm Riemann}(e_a)$ for
all $f\in \CA$, i.e. we need only give it on the basic 1-forms.
\eproof

\section{Lifting and Ricci curvature}

We now want to compute the Ricci tensor of our canonical spin
connection. As explained in \cite{Ma:rieq}, we need for this a
suitable lift $i:\Omega^2\to \Omega^1\tens_{\CA}\Omega^1$. We will
look for this as induced by a map $i:\Lambda^2\to
\Lambda^1\tens\Lambda^1$. For quantum groups we explain a natural
choice for the map $i$ and we compute it.

First of all, note that the metric $\eta$ defines a Hodge *
operator as in \cite{GomMa:coh}. The formulae can be obtained in
our conventions by a conversion process as in Section~2. Using
$*:\Lambda^2\to\Lambda^2$ we then decompose the invariant 2-forms
into self-dual and antiself-dual pieces
$\Lambda^2=\Lambda^2_+\oplus\Lambda^2_-$, where \eqn{selfdual}{
\Lambda^2_+=\{e_{ac},\ e_{bd}-\mu e_{ab},\ e_{ad}-q^{-2}e_{cb}\},
\quad \Lambda^2_-=\{e_{ab},\ e_{cd},\ e_{ad}+e_{cb}\}} and where
we use the shorthand $e_{ac}=e_a\wedge e_c$, etc. as in the proof
of Theorem~3.1. The $\Lambda^2_\pm$ each form 3-dimensional
comodules under $\Delta_R$. On the other hand the eigenvalues of
$\Psi$ in our case are
\[ {\rm Eigenvalues}(\Psi)=\{1,-q^{-2},-q^2\}\]
The 1-eigenvalue subspace is the kernel of
$\wedge:\Lambda^1\tens\Lambda^1\to \Lambda^2$. The
$-q^{\pm2}$-eigenvalue subspaces $E_\pm$ are again 3-dimensional
and are comodules under $\Delta_R$ since $\Psi$ is covariant.

\begin{propos} The map $i$ defined by the eigenspace decomposition of
$\Lambda^1\tens\Lambda^1$ under $\Psi$ is such that
$i(\Lambda^2_\pm)=E_\pm$, with
\begin{eqnarray*}
  i(e_{ac})&=&[2]_{q^2}^{-1}( q^2 e_a\tens e_c-e_c\tens e_a)\\
  i(e_{bd}-\mu e_{ab})&=&[2]_{q^2}^{-1} q^2(e_b\tens e_d
  -e_d\tens e_b-\mu e_a\tens e_b)\\
 i(e_{ad}-q^{-2}e_{cb})&=&[2]_{q^2}^{-1}( q^2 e_a\tens e_d
 - e_d\tens e_a+e_b\tens e_c-e_c\tens e_b-\mu e_a\tens e_a)\\
i(e_{ab})&=&[2]_{q^2}^{-1}(  e_a\tens e_b-q^2 e_b\tens e_a)\\
 i(e_{cd})&=&[2]_{q^2}^{-1}(  e_c\tens e_d-e_d\tens e_c
 +\mu e_c\tens e_a)\\
  i(e_{ad}+e_{cb})&=&[2]_{q^2}^{-1}( e_a\tens e_d-q^2 e_d\tens e_a
  +e_c\tens e_b-e_b\tens e_c+\mu e_a\tens e_a)
\end{eqnarray*}
\end{propos}
\proof The procedure to construct $i$ works for any quantum group
bicovariant calculus provided only that we can fully diagonalise
$\Psi:\Lambda^1\tens\Lambda^1\to \Lambda^1\tens\Lambda^1$. In this
case we let $M$ be the space spanned by all eigenvectors with
eigenvalues other than $1$.  Then the restriction $\wedge|_M:M\to
\Lambda^2$ is a linear isomorphism and we define
$i=(\wedge|_M)^{-1}$. We then compute in our example and find that
it maps (anti)-self dual 2-forms to the two eigenspaces as stated,
thereby giving these $-q^{\pm2}$-eigenspaces a geometrical
interpretation in terms of the Hodge $*$ operator. \eproof

Let us note that for low odd roots of unity we find that the
moduli space of all possible equivariant splittings
$i:\Lambda^2\to\Lambda^1\tens\Lambda^1$ of the wedge product is
2-dimensional. The above construction gives one point in this
moduli space but it works also for generic $q$ and is nonsingular
as $q\to 1$.

Next applying $i$ to the first factor of Riemann and for our form
of Riemann tensor where there is no functional dependence, we have
$i({\rm Riemann}(e_a)))\in \Lambda^1\tens
\Lambda^1\tens\Lambda^1$. We then define the Ricci tensor by a
`quantum trace' \eqn{ricci}{ {\rm Ricci}=\sum_{i,j,k}\<f^k,i({\rm
Riemann}(e_j)\>\eta^{ij}\eta_{ik}\in \Lambda^1\tens\Lambda^1}
where $\{f^i\}$ is a dual basis to $\{e_i\}$ and is evaluated
against the left most copy. Here $\eta=\eta^{ij}e_i\tens e_j$
(summation understood) defines the metric components for the
metric in Section~3 and $\eta_{ij}$ is the inverse matrix. Because
the matrix is not symmetric the combination $u^j{}_k\equiv
\eta^{ij}\eta_{ik}$ is not the identity matrix as it would be in a
usual trace. Instead we use a kind of `quantum trace' in defining
the Ricci. This is a small change from the formula in
\cite{Ma:rieq} (where we used the usual trace) and has the merit
that all constructions remain covariant under the background
quantum group (not only under a classical change of basis as
discussed there). The physical meaning of this improved Ricci
tensor is that we use the metric to raise the input index of
Riemann and then contract as usual with another (but transposed)
copy of the metric.

\begin{propos} The Ricci tensor for the canonical spin connection
in Theorem~3.1 is
\[ {\rm Ricci}=-{2q^2\over [4]_{q^2}}(\eta+{q^4\over
(1+q^2)}\theta\tens\theta)\]
i.e. an `Einstein space' up to a shift by $\theta\tens\theta$.
\end{propos}
\proof We first compute
\begin{eqnarray*} i(e_b\wedge e_z)&=&[2]_{q^2}^{-1}q^2((\mu-q^{-2})
e_a\tens e_b+e_b\tens e_a+e_d\tens e_b-e_b\tens e_d)\\
 i(e_c\wedge e_z)&=&[2]_{q^2}^{-1}((q^{-2}-\mu))e_c\tens e_a
 -e_a\tens e_c+e_d\tens e_c-e_c\tens e_d)\\
 i(e_c\wedge e_b)&=&[2]_{q^2}^{-2}q^22(e_c\tens e_b-e_b\tens
 e_c-{q^2\mu\over 2}(e_a\tens e_d+e_d\tens e_a)+\mu e_a\tens e_a).
 \end{eqnarray*}
 We also compute the quantum trace matrix as
\[ u=\begin{pmatrix}1&0&0&0\\ 0& q^2 & 0& 0\\ 0& 0& q^{-2} & 0\\
0&0&0&1\end{pmatrix}
\]
in the basis $\{e_a,e_c,e_b,e_d\}$. We then contract against the
Riemann tensor from Section~3. \eproof

Thus we find the same phenomenon as for the $S_3$ finite symmetric
group example in \cite{Ma:rieq}, i.e. the natural `Killing metric'
on $\C_q[SL_2]$ has covariantly constant Ricci tensor essentially
proportional to the metric up to a shift by $\theta\tens\theta$.
It is interesting that, as for $S_3$, this shift is exactly such
as to make the shifted metric degenerate. But for all nearby
metric we could take the shifted metric as a new metric without
changing most of the geometry. Also as in $S_3$ our contraction
conventions  (dictated by the noncommutative geometry) are such
that negative Ricci curvature corresponds to what in usual
geometry would be positive Ricci curvature. In this sense the
above result fits our expected picture of $\C_q[SU_2]$ with real
$q$ and a suitable $*$-structure as a noncommutative sphere.

\section{Dirac operator}

The additional ingredient for a Dirac operator is a choice of
spinor representation $W$ of the frame quantum group and
equivariant gamma-matrices $\gamma:\Lambda^1\to {\rm End}(W)$. The
spinor bundle in our case is just the tensor product $\CA\tens W$,
which is the space of spinors. We take the 2-dimensional
representation (i.e. a Weyl spinor) so a spinor has components
$\psi^\alpha\in\C_q[SL_2]$ for $\alpha=1,2$.

Since $\Lambda^1$ for the differential calculus of $\C_q[SL_2]$
was originally given in the endomorphism basis
$\{e_\alpha{}^\beta\}$, the canonical gamma-matrices proposed in
\cite{Ma:rieq} are just the identity map in that basis. Or in
terms of our above $\{e_i\}$ they provide the conversion according
to \eqn{gamma}{ \gamma(e_i)^\alpha{}_\beta e_\alpha{}^\beta
=e_i,\quad \gamma(e_a)=\begin{pmatrix}1&0\\0&0\end{pmatrix},\quad
\gamma(e_b)=\begin{pmatrix}0&1\\0&0\end{pmatrix},\quad {\rm etc.}}
If we take more usual linear combinations $e_x,e_y,e_z,\theta$
(where $e_x,e_y$ are linear combinations of $e_b,e_c$) then the
gamma-matrices would have a more usual form of Pauli matrices and
the identity, but this this is not particularly natural  when
$q\ne 1$ given that our metric is not symmetric.

Using the $\{e_\alpha{}^\beta\}$ basis, the canonical Dirac
operator in \cite{Ma:rieq} is \eqn{dirac}{
(\Dsl\psi)^\alpha=\del^\alpha{}_\beta\psi^\beta-A(\tilde\pi S^{-1}
t^\gamma{}_\beta)^\alpha{}_\gamma \psi^\beta} where
\[ A=A^\alpha{}_\beta e_\alpha{}^\beta,\quad \extd
f=\del^\alpha{}_\beta(f) e_\alpha{}^\beta,\quad \forall
A\in\Omega^1,\quad f\in \C_q[SL_2].\]

\begin{propos} The Dirac operator for the canonical
spin connection on $\C_q[SL_2]$ in Theorem~3.1 is
\[ \Dsl=\begin{pmatrix}\del^a+ {[3]\over [4]}& \del^b \\
\del^c& \del^d+ {[3]\over [4]} \end{pmatrix}=\dsl+{[3]\over [4]}\]
where $\dsl=\del^i\gamma(e_i)$ and $[n]=(q^n-q^{-n})/(q-q^{-1})$.
\end{propos}
\proof From (\ref{tildepi})  for $\tilde\pi$, and the specific
form of our spin connection from Theorem~3.1, we have
\begin{eqnarray*} && A(\tilde\pi S^{-1}a)=A(\tilde\pi
d)=-q^{-1}\alpha e_z,\quad
A(\tilde\pi S^{-1}b)=-q^{-1}A(\tilde\pi b)=-q^{-1}\beta e_c\\
&& A(\tilde\pi S^{-1}c)=-q A(\tilde\pi c)=-q\beta e_b,\quad
A(\tilde\pi S^{-1}d)=A(\tilde\pi a)=q\alpha e_z.\end{eqnarray*}
where $\alpha=q^4/[4]_{q^2}$ and $\beta=1/(q^2+q^{-2})$ are the
coefficients in Theorem~3.1. We then convert to the spinor basis
with $(e_i)^\alpha{}_\beta=\gamma(e_i){}^\alpha{}_\beta$ so, e.g.
$(e_b){}^1{}_2=1$ and its other components are zero. We then
compute the matrix $\Asl{}^\alpha{}_\beta=A(\tilde\pi S^{-1}
t^\gamma{}_\beta)^\alpha{}_\gamma$ so that $\Dsl=\dsl-\Asl$. We
have
\[  A(\tilde\pi S^{-1}t^1{}_1)^1{}_1+A(\tilde\pi
S^{-1}t^2{}_1)^1{}_2=-\alpha q^{-3}-q\beta,\quad A(\tilde\pi
S^{-1}t^1{}_2)^2{}_1+A(\tilde\pi
S^{-1}t^2{}_2)^2{}_2=-q^{-1}\beta-q\alpha\] and zero for the
off-diagonal entries of $\Asl$. Remarkably, these two diagonal
entries coincide and equal $-q[3]_{q^2}/[4]_{q^2}=-[3]/[4]$ on
inserting the values of $\alpha,\beta$. \eproof

We see that the Dirac operator is the naive $\dsl$ that one might
write guess without a spin connection, plus an offset of $[3]/[4]$
reflecting the constant background curvature. This is much the
same phenomenon as for $S_3$ in \cite{Ma:rieq} and in keeping with
the broad geometrical picture.

\begin{propos} At least for $r=3,5,7$ the above Dirac operator at
an $r$'th root of unity has $r$ eigenvalues of $\dsl$ lying
equally spaced on a circle in the complex plane. Here
\[ {\rm Eigenvalues}(\dsl)=\{{q^2[m]_q\over [2]_q},\ m=0,1,\cdots,r-1\}\]
The offset value in $\Dsl$ due to the constant curvature
background is
\[ {\sin({2\pi 3\over r})\over\sin({2\pi 4\over r})}=0,\
\sqrt{{\sqrt{5}-1\over\sqrt{5}+1}},\ -1,\]
for $r=3,5,7$. For $r=9$ it jumps to approximately 2.53... and
thereafter decays asymptotically towards ${3\over 4}$ as $r\to
\infty$.
\end{propos}
\proof  It is easiest (and most natural mathematically) to compute
the unnormalised $\dsl$ where $\extd=[\theta,\ \}$ without the
factor $\mu^{-1}$ that was inserted for the classical limit. The
eigenvalues in this case are computed using the explicit formulae
for $\extd$ in Section~2 and the corresponding $\del^a$ without
this factor, and one obtains $\{q^{m}-1|\ m=0,1,\cdots ,r-1\}$
i.e. all the conjugate roots of unity, shifted by -1. The
normalised $\dsl$ is then as stated. Note that we have one
eigenvalue for each irreducible representation of the reduced
quantum group $u_q(sl_2)$. \eproof

We expect the above results for all odd roots of unity and see
that the limit as $r\to \infty$ fits with the real $q\to 1$ limit.
The two real eigenvalues of $\dsl$ are 0,-1 (for $m=0,r-2$) and
the offset due to the curvature for large $r$ ensures that these
appear in $\Dsl$ on either side of zero. For small $r$ it is
interesting that the offset  vanishes for $r=3$, while for $r=5$
it is
\[ 4\cos^2({\pi\over 5})-2\]
which is reminiscent of the discrete series $\{4\cos^2({\pi\over
r})\}$ for the values of the Jones index of a type $II_1$
subfactor\cite{JonSun}. Likewise the distance between the $m=0$
and $m=1$ eigenvalues is $1/(2\cos({\pi\over r}))$, while the
radius of the circle of eigenvalues is a similar series
$|\mu|^{-1}=1/(2\sin({2\pi\over r}))$. Note, however, that these
observations, as well as the reality of the offset due to the
curvature, depend on the normalisation $\mu^{-1}$ in our initial
definition of $\extd$, which is not the only possible one. We note
also that Proposition~5.2 relating the eigenvalues of the Dirac
operator to the $q$-integer dimensions of irreducible
representations should be expected to apply also for the full
$\C_q[SL_2]$ at generic $q$, but with $m$ no longer truncated at
$r-1$. For generic $q$ a different approach to Dirac operators
appeared shortly after \cite{Ma:rieq} in \cite{Heck:spi}, with
different results than ours.

Finally, at roots of unity, the Dirac operator is {\em not} fully
diagonalisable, just as the Laplace operator was not in
\cite{GomMa:coh}. We demonstrate this for $r=3$:

\begin{propos} For $r=3$, the 14 zero-modes of $\dsl=\Dsl$ are
\[ \begin{pmatrix}0\\ 1\end{pmatrix},\  \begin{pmatrix}0\\ b^2d\end{pmatrix},
\ \begin{pmatrix}0\\ bd^2\end{pmatrix},\ \begin{pmatrix}-b\\
q a\end{pmatrix},\  \begin{pmatrix}-b^2d\\ b\end{pmatrix}, \
\begin{pmatrix}bd^2\\ d\end{pmatrix},\quad \begin{pmatrix}cd^2\\ qc^2d\end{pmatrix}\]
\[ \begin{pmatrix}1\\ 0\end{pmatrix},\  \begin{pmatrix}c^2a\\ 0\end{pmatrix},
\ \begin{pmatrix}ca^2\\ 0\end{pmatrix},\ \begin{pmatrix}q^2d\\
-c\end{pmatrix},\  \begin{pmatrix}c\\ -c^2a\end{pmatrix}, \
\begin{pmatrix}a\\ ca^2\end{pmatrix},\quad \begin{pmatrix}q^2b^2a\\ ba^2\end{pmatrix}\]
The 12 massive modes with eigenvalue $-1$ are:
\[ \begin{pmatrix}0\\ b\end{pmatrix},\  \begin{pmatrix}0\\ d\end{pmatrix},
\ \begin{pmatrix}0\\ b^2a\end{pmatrix},\ \begin{pmatrix}b\\
q^2a\end{pmatrix},\  \begin{pmatrix}1\\ -cd^2\end{pmatrix}, \
\begin{pmatrix}cd^2\\ -q^2c^2d\end{pmatrix}\]
\[ \begin{pmatrix}c\\ 0\end{pmatrix},\  \begin{pmatrix}a\\ 0\end{pmatrix},
\ \begin{pmatrix}c^2d\\ 0\end{pmatrix},\ \begin{pmatrix}qd\\
c\end{pmatrix},\  \begin{pmatrix}-ba^2\\ 1\end{pmatrix}, \
\begin{pmatrix}-qb^2a\\ ba^2\end{pmatrix}\]
The 12 massive modes with eigenvalue $q^2$ are:
\[ \begin{pmatrix}0\\ b^2\end{pmatrix},\  \begin{pmatrix}0\\ d^2\end{pmatrix},
\ \begin{pmatrix}0\\ bd\end{pmatrix},\ \begin{pmatrix}-b^2\\
qba\end{pmatrix},\  \begin{pmatrix}d^2\\ -q^2cd\end{pmatrix}, \
\begin{pmatrix}-q^2bd\\ cb-1\end{pmatrix}\]
\[ \begin{pmatrix}c^2\\ 0\end{pmatrix},\  \begin{pmatrix}a^2\\ 0\end{pmatrix},
\ \begin{pmatrix}ca\\ 0\end{pmatrix},\ \begin{pmatrix}q^2cd\\
-c^2\end{pmatrix},\  \begin{pmatrix}-qba\\ a^2\end{pmatrix}, \
\begin{pmatrix}cb-1\\ -qca\end{pmatrix}\]
\end{propos}
\proof We first compute the dimensions of the eigenspaces of
$\dsl$ as a $54\times 54$ matrix. Due to its size, we then treat
it as a numerical matrix to find approximate eigenvectors, which
we fit to powers of $q$ to obtain the required linearly
independent number of exact eigenvectors verified directly. We are
left with 16 non diagonalisable modes. We work in the
$\{c^kb^nd^m\}$ basis but simplify final results using the
identities $a=(1+q^2cb)d^2$ and $a^2=(1-qcb+q^2c^2b^2)d$. These
imply in particular that $c^2a^2=c^2d$, $c^2d^2=c^2a$,
$b^2a^2=b^2d$, $b^2d^2=b^2a$. \eproof

We see that the components of the modes of $\dsl$ with real
eigenvalues $0,-1$ are given by entries which are the 13 massless
spin zero modes $1, a, b, c, d, c^2d, cd^2, c^2a, ca^2, b^2d,
bd^2, b^2a, ba^2$ found in \cite{GomMa:coh}. This is a remnant of
the idea that the Dirac operator is the square root of the Laplace
one. On the other hand the components of the modes of $\dsl$ of
complex mass $q^2$ have entries from the massive spin zero modes
$a^2, b^2, c^2, d^2, cd, bd, ca, ba, cb-1$ in \cite{GomMa:coh}.
Note that at a mathematical level without the $\mu^{-1}$
normalisation in $\dsl$, the massive modes would have eigenvalues
$q-1,q^2-1$ and we see that these each have an equal number of
modes. Their physical make up, however, seems to be quite
different as just explained. Finally, we have written out the
modes in such a way as to manifest a different symmetry, namely
for each mass, they fall into pairs given by swapping the entries
of the matrix and replacing $a\leftrightarrow d$,
$b\leftrightarrow c$, $q\leftrightarrow q^{-1}$. This appears to
be a reasonable charge conjugation operation.

\bigskip

\baselineskip 14pt

\end{document}